# EIGENVALUES AND EIGENFUNCTIONS OF Q-DIRAC SYSTEM

FATMA HIRA[1]

***Abstract.*** *In this paper, we deal with a $q$-Dirac system. We investigate some spectral properties and the asymptotic behavior of the eigenvalues and the eigenfunctions of this $q$-Dirac system.*

## 1. INTRODUCTION

We consider a $q$-Dirac system which consists of the system of $q$-differential equations

$$\begin{cases} -\dfrac{1}{q} D_{q^{-1}} y_2(x) + p(x) y_1(x) = \lambda y_1(x), \\ D_q y_1(x) + r(x) y_2(x) = \lambda y_2(x), \end{cases} \quad (1)$$

and the boundary conditions

$$B_1(y) := k_{11} y_1(0) + k_{12} y_2(0) = 0, \quad (2)$$

$$B_2(y) := k_{21} y_1(a) + k_{22} y_2(aq^{-1}) = 0, \quad (3)$$

where $k_{ij} (i, j = 1, 2)$ are real numbers, $\lambda$ is a complex eigenvalue parameter, $y(x) = \begin{pmatrix} y_1(x) \\ y_2(x) \end{pmatrix}$, $0 \le x \le a < \infty$, $p(x)$ and $r(x)$ are real-valued functions defined on $[0, a]$ and continuous at zero and $p(x), r(x) \in L_q^1(0, a)$ (see[1]).

In[1], the authors introduced a $q$-analog of one-dimensonal Dirac equation (1) and they investigated the existence and uniqueness of the solution of this equation and also gave some spectral properties of the problem (1)-(3). Dissipative, accumulative, self-adjoint for similar $q$-Dirac equation were described in [2].

In this paper, we study similar spectral properties and obtain the asymptotic formulas of the eigenvalues and the eigenfunctions of the problem (1)-(3) in the light of the theory of $q$-(basic) Sturm-Liouville problems (see [3,4]).

## 2. PRELIMINARIES

In this section we introduce some of the required q-notations and results. Throughout this paper $q$ is a positive number with $0 < q < 1$.

A set $A \subseteq \mathbb{R}$ is called $q$-geometric if, for every $x \in A, qx \in A$. Let $f$ be a real or complex-valued function defined on a $q$-geometric set $A$. The $q$-difference operator is defined by

$$D_q f(x) := \frac{f(x) - f(qx)}{x(1-q)}, x \neq 0. \tag{4}$$

---

[1] Hitit University, Arts and Science Faculty, Department of Mathematics, 19030, Çorum, Turkey.
E-mail:fatmahira@yahoo.com.tr, fatmahira@hitit.edu.tr

If $0 \in A$, the $q$-derivative at zero is defined to be

$$D_q f(0) := \lim_{n \to \infty} \frac{f(xq^n) - f(0)}{xq^n}, x \in A, \tag{5}$$

if the limit exists and does not depend on $x$. Also, for $x \in A$, $D_{q^{-1}}$ is defined to be

$$D_{q^{-1}} f(x) := \begin{cases} \frac{f(x) - f(q^{-1}x)}{x(1-q^{-1})}, & x \in A \setminus \{0\}, \\ D_q f(0), & x = 0, \end{cases} \tag{6}$$

provided that $D_q f(0)$ exists. The following relation can be verified directly from the definition

$$D_{q^{-1}} y(x) = (D_q y)(xq^{-1}). \tag{7}$$

A right inverse, $q$-integration, of the $q$-difference operator $D_q$ is defined by Jackson [5] as

$$\int_0^x f(t) d_q t := x(1-q) \sum_{n=0}^{\infty} q^n f(xq^n), x \in A, \tag{8}$$

provided that the series converges. A $q$-analog of the fundamental theorem of calculus is given by

$$D_q \int_0^x f(t) d_q t = f(x), \quad \int_0^x D_q f(t) d_q t = f(x) - \lim_{n \to \infty} f(xq^n), \tag{9}$$

where $\lim_{n \to \infty} f(xq^n)$ can be replaced by $f(0)$ if $f$ is $q$-regular at zero, that is, if $\lim_{n \to \infty} f(xq^n) = f(0)$, for all $x \in A$. Throughout this paper, we deal only with functions $q$-regular at zero.

The $q$-type product formula is given by

$$D_q(fg)(x) = g(x) D_q f(x) + f(qx) D_q g(x), \tag{10}$$

and hence the $q$-integration by parts is given by

$$\int_0^a g(x) D_q f(x) d_q x = (fg)(a) - (fg)(0) - \int_0^a D_q g(x) f(qx) d_q x, \tag{11}$$

where $f$ and $g$ are $q$-regular at zero.

For more results and properties in $q$-calculus, readers are referred to the recent works [6-9].

The basic trigonometric functions $\cos(z;q)$ and $\sin(z;q)$ are defined on $\mathbb{C}$ by

$$\cos(z;q) := \sum_{n=0}^{\infty} \frac{(-1)^n q^{n^2} (z(1-q))^{2n}}{(q;q)_{2n}}, \tag{12}$$

$$\sin(z;q) := \sum_{n=0}^{\infty} \frac{(-1)^n q^{n(n+1)} (z(1-q))^{2n+1}}{(q;q)_{2n+1}}, \tag{13}$$

and they are $q$-analogs of the cosine and sine functions (see [7, 10, 11]).

**Theorem 2.1.** ([12]) If $\{x_m\}$ and $\{y_m\}$ are the positive zeros of $\cos(z;q)$ and $\sin(z;q)$, respectively, then we have for sufficiently large m,

$$\{x_m\} = q^{-m+1/2}(1-q)^{-1}(1+O(q^m)), \tag{14}$$

$$\{y_m\} = q^{-m}(1-q)^{-1}(1+O(q^m)). \tag{15}$$

**Corollary 2.1.** ([12, Corollaries 3.2 and 3.4]) For $r := |z| \to \infty$ we have

$$M(r; \cos(z;q)) = O\left(\exp\left(-\frac{(\log r(1-q))^2}{\log q}\right)\right), \tag{16}$$

$$M(r; \sin(z;q)) = O\left(\exp\left(-\frac{(\log r(1-q))^2}{\log q}\right)\right). \tag{17}$$

## 3. SPECTRAL PROPERTIES AND ASYMPTOTIC FORMULAS

In this section we give some spectral properties similar to [1,13], then we obtain asymptotic formulas for the eigenvalues and the eigenfunctions of the problem (1)-(3).

**Lemma 3.1.** The eigenfunctions $y(x, \lambda_1)$ and $z(x, \lambda_2)$ corresponding to different eigenvalues $\lambda_1 \neq \lambda_2$ are orthogonal, i.e.,

$$\int_0^a y^\perp z \, d_q x = \int_0^a \{y_1(x,\lambda_1) z_1(x,\lambda_2) + y_2(x,\lambda_1) z_2(x,\lambda_2)\} d_q x = 0, \tag{18}$$

where $y^\perp = (y_1, y_2)$.

*Proof.* Since $y(x, \lambda_1)$ and $z(x, \lambda_2)$ are solutions of the $q$-system (1),

$$-\frac{1}{q} D_{q^{-1}} y_2 + \{p(x) - \lambda_1\} y_1 = 0,$$

$$D_q y_1 + \{r(x) - \lambda_1\} y_2 = 0,$$

$$-\frac{1}{q} D_{q^{-1}} z_2 + \{p(x) - \lambda_2\} z_1 = 0,$$

$$D_q z_1 + \{r(x) - \lambda_2\} z_2 = 0.$$

Multiplying by $z_1, z_2, -y_1$ and $-y_2$, respectively, and adding together and also using the formulas (7) and (9), we obtain

$$D_q\left(y_1(x,\lambda_1) z_2(xq^{-1}, \lambda_2) - y_2(xq^{-1}, \lambda_1) z_1(x, \lambda_2)\right)$$
$$= (\lambda_1 - \lambda_2)\{y_1(x,\lambda_1) z_1(x,\lambda_2) + y_2(x,\lambda_1) z_2(x,\lambda_2)\}. \tag{19}$$

Then integrating from $0$ to $a$, we have

$$(\lambda_1 - \lambda_2)\int_0^a y^\perp(x,\lambda_1) z(x,\lambda_2) d_q x = \left\{y_1(x,\lambda_1) z_2(xq^{-1}, \lambda_2) - y_2(xq^{-1}, \lambda_1) z_1(x, \lambda_2)\right\}\Big|_0^a. \tag{20}$$

It follows from the boundary conditions (2) and (3) the right hand side vanishes. Therefore, we get

$$(\lambda_1 - \lambda_2)\int_0^a y^\perp(x,\lambda_1) z(x,\lambda_2) d_q x = 0. \tag{21}$$

The lemma is thus proved, since $\lambda_1 \neq \lambda_2$.    □

**Lemma 3.2.** The eigenvalues of the problem (1)-(3) are real.

*Proof.* Assume the contrary that $\lambda_0$ is a nonreal eigenvalues of the problem (1)-(3). Let $y(x,\lambda_0)$ be a corresponding (nontrivial) eigenfunctions. $\overline{\lambda_0}$ is also an eigenvalue, corresponding to the eigenfunctions $\overline{y}(x,\lambda_0)$. Since $\lambda_0 \neq \overline{\lambda_0}$ by the previous lemma,

$$\int_0^a \left( |y_1(x,\lambda_0)|^2 + |y_2(x,\lambda_0)|^2 \right) d_q x = 0, \tag{22}$$

hence, $y(x,\lambda_0) \equiv 0$ and this is a contradiction. Consequently, $\lambda_0$ must be real. □

Now, we will construct a special system of solution of $q$-system (1). Let $\phi(x,\lambda) = \begin{pmatrix} \phi_1(.,\lambda) \\ \phi_2(.,\lambda) \end{pmatrix}$ be a solution of the $q$-system (1) that satisfies the initial conditions

$$\phi_1(0,\lambda) = k_{12}, \quad \phi_2(0,\lambda) = -k_{11}. \tag{23}$$

The existence and uniqueness of this solution for the $q$-system (1) were presented in [1]. It is obvious that $\phi(x,\lambda)$ satisfies the boundary condition (2).

**Theorem 3.1.** The following integral equations hold for the solution of $\phi(x,\lambda)$

$$\phi_1(x,\lambda) = k_{12}\cos(\lambda x;q) - k_{11}\sin(\lambda x;q)$$
$$+ q\int_0^x \{\sin(\lambda x;q)\cos(\lambda qt;q) - \cos(\lambda x;q)\sin(\lambda qt;q)\} p(qt)\phi_1(qt,\lambda)d_q t \tag{24}$$
$$- \int_0^x \{\cos(\lambda x;q)\cos(\lambda\sqrt{q}t;q) + \sqrt{q}\sin(\lambda x;q)\sin(\lambda\sqrt{q}t;q)\} r(t)\phi_2(t,\lambda)d_q t,$$

$$\phi_2(x,\lambda) = -k_{12}\sqrt{q}\sin(\lambda\sqrt{q}x;q) - k_{11}\cos(\lambda\sqrt{q}x;q)$$
$$+ q\int_0^x \{\cos(\lambda\sqrt{q}x;q)\cos(\lambda qt;q) + \sqrt{q}\sin(\lambda\sqrt{q}x;q)\sin(\lambda qt;q)\} p(qt)\phi_1(qt,\lambda)d_q t \tag{25}$$
$$+ \sqrt{q}\int_0^x \{\sin(\lambda\sqrt{q}x;q)\cos(\lambda\sqrt{q}t;q) - \cos(\lambda\sqrt{q}x;q)\sin(\lambda\sqrt{q}t;q)\} r(t)\phi_2(t,\lambda)d_q t.$$

*Proof.* Let us construct two solutions of $q$-system (1) as

$$\varphi_1(.,\lambda) = \begin{pmatrix} \varphi_{11}(x,\lambda) \\ \varphi_{12}(x,\lambda) \end{pmatrix} = \begin{pmatrix} \cos(\lambda x;q) \\ -\sqrt{q}\sin(\lambda\sqrt{q}x;q) \end{pmatrix},$$
$$\varphi_2(.,\lambda) = \begin{pmatrix} \varphi_{21}(x,\lambda) \\ \varphi_{22}(x,\lambda) \end{pmatrix} = \begin{pmatrix} \sin(\lambda x;q) \\ \cos(\lambda\sqrt{q}x;q) \end{pmatrix}, \tag{26}$$

for $p(x) = r(x) \equiv 0$, with the Wronskian

$$W(\varphi_1,\varphi_2)(x,\lambda) = \varphi_{11}(x,\lambda)\varphi_{22}(xq^{-1},\lambda) - \varphi_{21}(x,\lambda)\varphi_{12}(xq^{-1},\lambda) = 1. \tag{27}$$

The function

$$y(.,\lambda) = \begin{pmatrix} y_1(x,\lambda) \\ y_2(x,\lambda) \end{pmatrix} = \begin{pmatrix} c_1\cos(\lambda x;q) + c_2\sin(\lambda x;q) \\ -c_1\sqrt{q}\sin(\lambda\sqrt{q}x;q) + c_2\cos(\lambda\sqrt{q}x;q) \end{pmatrix} \tag{28}$$

is a fundamental set of the $q$-system (1) for $p(x) = r(x) \equiv 0$. Using $q$-analogue of the method of variation of constants, a particular solution of the $q$-system (1) may be given by

$$\begin{cases} y_1(x,\lambda) = c_1(x)\cos(\lambda x; q) + c_2(x)\sin(\lambda x; q), \\ y_2(x,\lambda) = -c_1(x)\sqrt{q}\sin(\lambda\sqrt{q}x; q) + c_2(x)\cos(\lambda\sqrt{q}x; q). \end{cases} \quad (29)$$

Hence the functions $c_i(x)(i=1,2)$ satisfy the $q$-linear system of equations

$$\begin{cases} \cos(\lambda x; q) D_{q^{-1}} c_1(x) + \sin(\lambda x; q) D_{q^{-1}} c_2(x) = -r(xq^{-1}) y_2(xq^{-1}, \lambda), \\ \sqrt{q}\sin(\lambda q^{-1/2} x; q) D_{q^{-1}} c_1(x) - \cos(\lambda q^{-1/2} x; q) D_{q^{-1}} c_2(x) = -qp(x) y_1(x, \lambda). \end{cases} \quad (30)$$

Since the equality (27) satisfies, (30) has a unique solution which leads

$$\begin{aligned} D_{q^{-1}} c_1(x) &= -r(xq^{-1})\cos(\lambda q^{-1/2}x; q) y_2(xq^{-1}, \lambda) - qp(x)\sin(\lambda x; q) y_1(x, \lambda), \\ D_{q^{-1}} c_2(x) &= qp(x)\cos(\lambda x; q) y_1(x, \lambda) - r(xq^{-1})\sqrt{q}\sin(\lambda q^{-1/2}x; q) y_2(xq^{-1}, \lambda). \end{aligned} \quad (31)$$

Using the formula (7) and replacing $x$ by $xq$ in (31), then we obtain

$$\begin{aligned} c_1(x) &= c_1 - q\int_0^x p(qt)\sin(\lambda qt; q) y_1(qt, \lambda) d_q t - \int_0^x r(t)\cos(\lambda\sqrt{q}t; q) y_2(t, \lambda) d_q t, \\ c_2(x) &= c_2 + q\int_0^x p(qt)\cos(\lambda qt; q) y_1(qt, \lambda) d_q t - \int_0^x r(t)\sqrt{q}\sin(\lambda\sqrt{q}t; q) y_2(t, \lambda) d_q t, \end{aligned} \quad (32)$$

when $c_i(x)(i=1,2)$ are $q$-regular at zero. That is the general solution $y(x,\lambda) = \begin{pmatrix} y_1(x,\lambda) \\ y_2(x,\lambda) \end{pmatrix}$ of $q$-system (1) is obtained to be

$$\begin{aligned} y_1(x,\lambda) &= c_1 \cos(\lambda x; q) + c_2 \sin(\lambda x; q) \\ &+ q\int_0^x \{\sin(\lambda x; q)\cos(\lambda qt; q) - \cos(\lambda x; q)\sin(\lambda qt; q)\} p(qt) y_1(qt, \lambda) d_q t \\ &- \int_0^x \{\cos(\lambda x; q)\cos(\lambda\sqrt{q}t; q) + \sqrt{q}\sin(\lambda x; q)\sin(\lambda\sqrt{q}t; q)\} r(t) y_2(t, \lambda) d_q t, \end{aligned} \quad (33)$$

$$\begin{aligned} y_2(x,\lambda) &= -c_1\sqrt{q}\sin(\lambda\sqrt{q}x; q) + c_2\cos(\lambda\sqrt{q}x; q) + \\ &+ q\int_0^x \{\cos(\lambda\sqrt{q}x; q)\cos(\lambda qt; q) + \sqrt{q}\sin(\lambda\sqrt{q}x; q)\sin(\lambda qt; q)\} p(qt) y_1(qt, \lambda) d_q t \\ &+ \sqrt{q}\int_0^x \{\sin(\lambda\sqrt{q}x; q)\cos(\lambda\sqrt{q}t; q) - \cos(\lambda\sqrt{q}x; q)\sin(\lambda\sqrt{q}t; q)\} r(t) y_2(t, \lambda) d_q t. \end{aligned} \quad (34)$$

It is easy to determine $c_1, c_2$ for which $\phi(x, \lambda)$ satisfies the $q$-system (1) and the conditions (23), then we obtain (24) and (25). □

**Theorem 3.2.** As $|\lambda| \to \infty$, the function $\phi(x, \lambda)$ has the following asymptotic relations

$$\phi_1(x, \lambda) = k_{12}\cos(\lambda x; q) - k_{11}\sin(\lambda x; q) + O\left(|\lambda|^{-1}\exp\left(\frac{-(\log|\lambda|x(1-q))^2}{\log q}\right)\right), \quad (35)$$

$$\phi_2(x,\lambda) = -k_{12}\sqrt{q}\sin(\lambda\sqrt{q}x;q) - k_{11}\cos(\lambda\sqrt{q}x;q)$$
$$+ O\left(|\lambda|^{-1}\exp\left(\frac{-(\log|\lambda|q^{1/2}x(1-q))^2}{\log q}\right)\right), \quad (36)$$

where for each $x \in (0,a]$ the $O$-terms are uniform on $\{xq^n : n \in \mathbb{N}\}$.

*Proof.* Similar to asymptotic relations for $q$-Sturm-Liouville problems in [4] and from Corollary (2.1), (35) and (36) can be obtained easily. □

**Lemma 3.3.** The eigenvalues of the problem (1)-(3) are simple.

*Proof.* The solution $\phi(x,\lambda)$ defined above is a nontrivial solution of $q$-system (1) satisfying the boundary condition (2). To find the eigenvalues of the problem (1)-(3), we have to insert this function into the boundary condition (3) and find the roots of the obtained equation. So, putting the function $\phi(x,\lambda)$ into the boundary condition (3) we get the following equation

$$\Delta(\lambda) = k_{21}\phi_1(a,\lambda) + k_{22}\phi_2(aq^{-1},\lambda). \quad (37)$$

Then $\dfrac{d\Delta(\lambda)}{d\lambda} = k_{21}\dfrac{\partial\phi_1(a,\lambda)}{\partial\lambda} + k_{22}\dfrac{\partial\phi_2(aq^{-1},\lambda)}{\partial\lambda}$. Let $\lambda_0$ be a double eigenvalue, and $\phi^0(x,\lambda_0)$ one of the corresponding eigenfunctions. Then the conditions $\Delta(\lambda_0) = 0$, $\dfrac{d\Delta(\lambda_0)}{d\lambda} = 0$ should be fulfilled simultaneously, i.e.,

$$k_{21}\phi_1^0(a,\lambda_0) + k_{22}\phi_2^0(aq^{-1},\lambda_0) = 0,$$
$$k_{21}\frac{\partial}{\partial\lambda}\phi_1^0(a,\lambda_0) + k_{22}\frac{\partial}{\partial\lambda}\phi_2^0(aq^{-1},\lambda_0) = 0. \quad (38)$$

Since $k_{21}$ and $k_{22}$ cannot vanish simultaneously, it follows from (38) that

$$\phi_1^0(a,\lambda_0)\frac{\partial\phi_2^0(aq^{-1},\lambda_0)}{\partial\lambda} - \phi_2^0(aq^{-1},\lambda_0)\frac{\partial\phi_1^0(a,\lambda_0)}{\partial\lambda} = 0. \quad (39)$$

Now, differentiating the $q$-system (1) with respect to $\lambda$, we obtain

$$-\frac{1}{q}D_{q^{-1}}\left(\frac{\partial y_2}{\partial\lambda}\right) + \{p(x) - \lambda\}\frac{\partial y_1}{\partial\lambda} = y_1,$$
$$D_q\left(\frac{\partial y_1}{\partial\lambda}\right) + \{r(x) - \lambda\}\frac{\partial y_2}{\partial\lambda} = y_2. \quad (40)$$

Multiplying the $q$-system (1) and (40) by $\dfrac{\partial y_1}{\partial\lambda}, \dfrac{\partial y_2}{\partial\lambda}, -y_1$ and $-y_2$, respectively, adding them together and integrating with respect to $x$ from $0$ and $a$, we obtain

$$\left\{y_2(xq^{-1},\lambda)\frac{\partial y_1(x,\lambda)}{\partial\lambda} - y_1(x,\lambda)\frac{\partial y_2(xq^{-1},\lambda)}{\partial\lambda}\right\}_0^a = \int_0^a \{y_1^2(x,\lambda) + y_2^2(x,\lambda)\}d_qx. \quad (41)$$

Putting $\lambda = \lambda_0$, taking into account that $\left.\dfrac{\partial\phi_1^0(x,\lambda_0)}{\partial\lambda}\right|_{x=0} = \left.\dfrac{\partial\phi_2^0(x,\lambda_0)}{\partial\lambda}\right|_{x=0} = 0$ by (24) and (25), and using the equality (39), we obtain the relation

$$\int_0^a \left\{\left(\phi_1^0(x,\lambda_0)\right)^2 + \left(\phi_2^0(x,\lambda_0)\right)^2\right\}d_qx = 0. \quad (42)$$

Hence $\phi_1^0(x,\lambda_0) = \phi_2^0(x,\lambda_0) \equiv 0$, which is impossible. Consequently $\lambda_0$ must be a simple eigenvalue. □

**Theorem 3.3.** As $|\lambda| \to \infty$ the function $\Delta(\lambda)$ has the following asymptotic relation

$$\Delta(\lambda) = k_{21}\left\{k_{12}\cos(\lambda a; q) - k_{11}\sin(\lambda a; q) + O\left(|\lambda|^{-1}\exp\left(\frac{-(\log|\lambda|a(1-q))^2}{\log q}\right)\right)\right\}$$

$$+ k_{22}\left\{-k_{12}\sqrt{q}\sin(\lambda q^{-1/2}a; q) - k_{11}\cos(\lambda q^{-1/2}a; q) + O\left(|\lambda|^{-1}\exp\left(\frac{-(\log|\lambda|aq^{-1/2}(1-q))^2}{\log q}\right)\right)\right\}.$$

(43)

*Proof.* The proof is immediate by substituting (35) and (36) into the relation

$$\Delta(\lambda) = k_{21}\phi_1(a,\lambda) + k_{22}\phi_2(aq^{-1},\lambda).$$
□

**Theorem 3.4.** The eigenvalues $\{\lambda_m\}$ are the zeros of $\Delta(\lambda)$ has the following asymptotic relations as $m \to \infty$:

*Case 1.* $k_{12} \neq 0, k_{11} = 0$;

$$i)\lambda_m = \frac{q^{-m+1/2}}{a(1-q)}\left(1 + O(q^{m/2})\right), \quad k_{21} = 0, \tag{44}$$

$$ii)\lambda_m = \frac{q^{-m+1/2}}{a(1-q)}\left(1 + O(q^m)\right), \quad k_{22} = 0, \tag{45}$$

*Case 2.* $k_{12} = 0, k_{11} \neq 0$;

$$i)\lambda_m = \frac{q^{-m+1}}{a(1-q)}\left(1 + O(q^{m/2})\right), \quad k_{21} = 0, \tag{46}$$

$$ii)\lambda_m = \frac{q^{-m}}{a(1-q)}\left(1 + O(q^m)\right), \quad k_{22} = 0. \tag{47}$$

*Proof.* Similar to asymptotic relation for $q$-Sturm-Liouville problems [4] and from Theorem 2.1 the asymptotic relations (44)-(47) can be obtained easily. □

Then from (35) and (36) and above theorem, the asymptotic relation of the eigenfunction of the problem (1)-(3) is given by

*Case 1.* $k_{12} \neq 0, k_{11} = 0$;

$$\phi(x,\lambda_m) = \begin{pmatrix} \phi_1(x,\lambda_m) \\ \phi_2(x,\lambda_m) \end{pmatrix}$$

$$= \begin{cases} k_{12}\cos(\lambda_m x; q) + O\left(|\lambda_m|^{-1}\exp\left(\frac{-(\log|\lambda_m|x(1-q))^2}{\log q}\right)\right), \\ -k_{12}\sqrt{q}\sin(\lambda_m\sqrt{q}x; q) + O\left(|\lambda_m|^{-1}\exp\left(\frac{-(\log|\lambda_m|q^{1/2}x(1-q))^2}{\log q}\right)\right), \end{cases}$$

*Case 2.* $k_{12} = 0, k_{11} \neq 0$;

$$\phi(x,\lambda_m) = \begin{pmatrix} \phi_1(x,\lambda_m) \\ \phi_2(x,\lambda_m) \end{pmatrix}$$

$$= \begin{cases} -k_{11}\sin(\lambda_m x;q) + O\left(|\lambda_m|^{-1}\exp\left(\frac{-(\log|\lambda_m|x(1-q))^2}{\log q}\right)\right), \\ -k_{11}\cos(\lambda_m\sqrt{q}x;q) + O\left(|\lambda_m|^{-1}\exp\left(\frac{-(\log|\lambda_m|q^{1/2}x(1-q))^2}{\log q}\right)\right). \end{cases}$$